\documentclass{article}

\usepackage{
 amsmath,
 amsxtra,
 amsthm,
 amssymb,
 etex,
 mathrsfs,
 }
\usepackage[all]{xy}

\newtheorem{theorem}{Theorem}[section]
\newtheorem{lemma}[theorem]{Lemma}

\newtheorem{proposition}[theorem]{Proposition}
\newtheorem{corollary}[theorem]{Corollary}

\theoremstyle{remark}
\newtheorem{remark}[theorem]{Remark}

\newtheorem{defn}[theorem]{Definition}

\newcommand{\bd}{\begin{defn}}
\newcommand{\ed}{\end{defn}}
\newcommand{\bl}{\begin{lemma}}
\newcommand{\el}{\end{lemma}}
\newcommand{\bp}{\begin{proposition}}
\newcommand{\ep}{\end{proposition}}
\newcommand{\bt}{\begin{theorem}}
\newcommand{\et}{\end{theorem}}
\newcommand{\bc}{\begin{corollary}}
\newcommand{\ec}{\end{corollary}}
\newcommand{\br}{\begin{remark}}
\newcommand{\er}{\end{remark}}
\newcommand{\ba}{\begin{array}}
\newcommand{\ea}{\end{array}}
\newcommand{\bpf}{\begin{proof}}
\newcommand{\epf}{\end{proof}}

\newcommand{\Q}{\mathbb{Q}}
\newcommand{\Qp}{\mathbb{Q}_p}
\newcommand{\Z}{\mathbb{Z}}
\newcommand{\Zp}{\mathbb{Z}_p}

\DeclareMathOperator{\Gal}{Gal}

\DeclareMathOperator{\coker}{coker}

\newcommand{\lra}{\longrightarrow}

\newcommand{\ot}{\otimes}

\newcommand{\ilim}{\displaystyle \mathop{\varinjlim}\limits}

\numberwithin{equation}{section}

\begin{document}
\title{Norm principle for even $K$-groups of number fields}
 \author{
  Meng Fai Lim\footnote{School of Mathematics and Statistics $\&$ Hubei Key Laboratory of Mathematical Sciences,
Central China Normal University, Wuhan, 430079, P.R.China.
 E-mail: \texttt{limmf@ccnu.edu.cn}} }
\date{}
\maketitle

\begin{abstract} \footnotesize
\noindent We investigate the norm maps of algebraic even $K$-groups of finite extensions of number fields. Namely, we show that they are surjective in most situations. In the event that they are not surjective, we give a criterion in determining when an element in the even $K$-group 
of the base field comes from a norm of an element from the even $K$-groups 
of the extension field. This latter criterion is only reliant on the real primes of the base field.


\end{abstract}

\section{Introduction}

The goal of the paper is study the norm map (also called the transfer map) of $K_{2n}$-groups of finite extension of number fields. The origin of this problem stems from the classical Hasse norm theorem which says that for a cyclic extension of number fields $L/F$, an element in $F^\times$ is a norm of an element of $L^\times$ if and only if it is a norm everywhere locally. This can be recast to a statement on $K_1$-groups of number fields, namely, a criterion on when an element in $K_1(F)$ can be realized as a norm of an element of $K_1(L)$. Motivated by this interpretation, Bak and Rehmann formulated and proved a $K_2$-analog of this statement in \cite{BR} (also see \cite{CT, MS}). One striking phenomenon was that the $K_2$-analog is valid for all finite extensions unlike its classical counterpart. Their results have then led them to the question of higher $K_m$-analogs. Although we will not address this in our paper, we should mention that had one replaced the algebraic $K$-groups by Milnor $K$-groups, such norm principles have been established in the works of Orlov-Vishik-Voevodsky \cite{OVV} and
Merkurjev-Suslin \cite{MS}. In \cite{Os}, {\O}stv{\ae}r considered the case of $K_{2n}$ for $n\geq 1$, where he succeeded in establishing these higher analog for the $2$-primary part of the $K_{2n}$-group. The goal of the present paper is to complete this study building on ideas of {\O}stv{\ae}r's in utilizing tools developed from recent progress on algebraic $K$-theory \cite{Vo}.

We begin introducing some preliminary notation to state the main results of this paper. For a ring $R$ with identity, write $K_m(R)$ for the algebraic $K$-groups of $R$ in the
sense of Quillen (see \cite{Qui72, Qui73a, Qui73b}). If $F$ is a number field, and $v$ is a (possibly archimedean) prime of $F$, we let $F_v$ denote the completion of $F$ at $v$. For a finite extension $L$ of $F$, denote by $N_{L/F}$ the norm map (also called the transfer map) $K_{2n}(L)\lra K_{2n}(F)$ on the $K$-groups (cf.\ \cite[Chap. IV, Definition 6.3.2]{WeiKbook}).

Let $v$ be a real prime of $F$. Denote by $l_v:F\lra F_v=\mathbb{R}$ the embedding at $v$. This in turn induces a homomorphism $K_{2n}(F)\lra K_{2n}(F_v)$ which by abuse of notation is also denoted as $l_v$. Since $K_{2n}(F)$ is a torsion group (for instance, see Lemma \ref{torsionK2}), the image of $l_v$ lies in $K_{2n}(F_v)_{tor}$, where we write $M_{tor}$ for the torsion subgroup of an abelian group $M$.

Our main result is as follows (compare with \cite[Theorem 2]{BR}).

\bt \label{main combine}
Let $L$ be a finite extension of $F$. Then we have the following exact sequence
\[ K_{2n}(L)\stackrel{N_{L/F}}{\lra} K_{2n}(F) \stackrel{l}{\lra} \prod_{v\in S_r}K_{2n}(F_v)_{tor}\lra 0,\]
where $S_r$ is the set of real primes of $F$ ramified in $L$ and the map $l$ is defined by $l(x) = (l_v(x))_{v\in S_r}$.
\et


Note that in the event that $n\neq 1~(\mathrm{mod}~4)$ or no real prime of $F$ is ramified in $L$, then Theorem \ref{main combine} is saying that $N_{L/F}$ is surjective. In other words, every element of $K_{2n}(F)$ is automatically a norm of an element in $K_{2n}(L)$ in these situations. For the remaining case, we have the following.

\bt \label{main cor}
Suppose that $n= 1~(mod~4)$ and that at least one real prime of $F$ is ramified in $L$. Let $x\in K_{2n}(F)$. Then $x$ is a norm of an element in $K_{2n}(L)$ if and only if $l_v(x) =0$ for every real prime of $F$ which is ramified in $L$.
\et

In particular, Theorem \ref{main cor} is saying that the property of $x$ being a norm of an element in $K_{2n}(L)$ is more or less reliant only on $F$. To illustrate this point, we record the following. Define $\mathcal{L}$ to be the collection of all finite extensions $L$ of $F$ such that $L$ has no real primes. Then the following is an immediate consequence of Theorem \ref{main cor}.

\bc \label{main corIllustrate}
Suppose that $n= 1~(mod~4)$ and that the number field $F$ has at least one real prime. Let $x\in K_{2n}(F)$. If there exists an $L_0\in \mathcal{L}$ such that $x$ is a norm of an element in $K_{2n}(L_0)$, then $x$ can be realized as a norm of an element in $K_{2n}(L)$ for every $L\in \mathcal{L}$. Vice versa, if there exists an $L_0\in \mathcal{L}$ such that $x$ is not a norm of an element in $K_{2n}(L_0)$, then $x$ cannot be realized as a norm element from $K_{2n}(L)$ for any $L\in \mathcal{L}$.
\ec

 Finally, in view of Theorem \ref{main cor}, it is natural to ask if one has an ``explicit'' way in determining when $l_v(x) =0$ for a real prime of $F$ that is ramified in $L$. For this, we provide such a possible criterion in the form of Propositions \ref{normK2} and \ref{normK2x} for the case $n=1$.

\section{Proof of the main results}

We shall prove our main results in this section. As a start, we have the following lemma.

\bl \label{torsionK2}
For $n\geq 1$, the group $K_{2n}(F)$ is a torsion group and has no nonzero $p$-divisible subgroups for all primes $p$.
\el

\bpf
This is certainly well-known. For the convenience of the readers, we sketch an  argument. By \cite[Theorem 3]{Sou}, one has a short exact sequence
\[ 0\lra K_{2n}(\mathcal{O}_F) \lra K_{2n}(F)\lra \bigoplus_v K_{2n-1}(k_v) \lra 0, \]
where $\mathcal{O}_F$ is the ring of integers of $F$ and the direct sum runs through all the finite primes of $F$ with $k_v$ being the residue field of $F_v$. The calculations of Quillen \cite{Qui73b} and Borel \cite{Bo} tell us that $K_{2n}(\mathcal{O}_F)$ is finite, where we note that the finiteness of $K_2$ was previously established by Garland \cite{Gar}. Quillen \cite{Qui72} has also established the finiteness of $K_{2n-1}(k_v)$. In view of these deep facts, we see that $K_{2n}(F)$ is an extension of a finite group and a sum of finite groups. The assertion of the lemma is now immediate from this.
\epf

We now prove the following which is essentially the $p$-primary part of our main results for odd $p$. For a given abelian group $N$, denote by $N[p^j]$ the subgroup of $N$ consisting of elements annihilated by $p^j$. One then writes $N[p^\infty]= \cup_{j\geq 1}N[p^j]$.

\bp \label{surjective main at odd p}
Suppose that $L/F$ is a finite extension of number fields. Then the
norm map
\[K_{2n}(L)[p^\infty]\lra K_{2n}(F)[p^\infty]\]
is surjective for every odd prime $p$. Furthermore, if  $L$ is a finite Galois extension of $F$, then the norm map induces an isomorphism
\[\big(K_{2n}(L)[p^\infty]\big)_{\Gal(L/F)}\cong K_{2n}(F)[p^\infty].\]
\ep

\bpf Observe that if $M$ is a finite Galois extension of $F$ containing $L$, it then follows from the transitivity of the norm maps that the surjectivity of $N_{L/F}$ is a consequence of the surjectivity of $N_{M/F}$. Hence we might as well assume that $L$ is a finite Galois extension of $F$ at the start, and to lighten notation, we shall write $G=\Gal(L/F)$. By the Universal Coefficient Theorem \cite[Chap.\ IV Theorem 2.5]{WeiKbook}, there is a short exact sequence
\[ 0\lra K_{2n+1}(F)/p^i \lra K_{2n+1}(F,\Z/p^i) \lra K_{2n}(F)[p^i]\lra 0,\]
where $K_{2n}(F,\Z/p^i)$ is the $K$-group with finite coefficients (for instance, see \cite{Br} or \cite[Chap.\ IV Section 2]{WeiKbook}). Taking direct limit, we obtain a short exact sequence
\[ 0\lra K_{2n+1}(F)\ot\Qp/\Zp \lra K_{2n+1}(F,\Qp/\Zp) \lra K_{2n}(F)[p^\infty]\lra 0.\]
One has a similar exact sequence for $L$, and together, they fit into the following commutative diagram
\[   \entrymodifiers={!! <0pt, .8ex>+} \SelectTips{eu}{}\xymatrix{
     &\big(K_{2n+1}(L)\ot\Qp/\Zp \big)_G  \ar[d] \ar[r] & K_{2n+1}(L,\Qp/\Zp)_G
    \ar[d] \ar[r] & \big(K_{2n}(L)[p^\infty]\big)_G \ar[d] \ar[r]& 0 \\
    0 \ar[r]^{} & K_{2n+1}(F)\ot\Qp/\Zp \ar[r]^{} & K_{2n+1}(F,\Qp/\Zp) \ar[r] &K_{2n}(F)[p^\infty] \ar[r]& 0  } \]
with exact rows and vertical maps induced by the norm maps. We claim that both assertions of the proposition will follow once we can show that the middle vertical map is an isomorphism. Indeed, supposing that for now the middle vertical map is an isomorphism. Then by a snake lemma argument, the rightmost vertical map is surjective with kernel being isomorphic to the cokernel of the leftmost vertical map. On the one hand, it follows from Lemma \ref{torsionK2} that this kernel is a torsion group with no nonzero $p$-divisible subgroups. On the other hand, being isomorphic to the cokernel of the leftmost vertical map, it has to be $p$-divisible. In conclusion, this forces the kernel to be trivial, and so the rightmost vertical map is an isomorphism yielding the conclusions of the proposition.

Therefore, it remains to verify that the middle vertical map is an isomorphism. For this, we need to recall some more terminologies. Denoting by $\mu_{p^n}$ the cyclic group generated by a primitive $p^n$-root of unity, we then write $\mu_{p^\infty}$ for the direct limit of the groups $\mu_{p^n}$. These have natural $\Gal(\bar{F}/F)$-module structures. Furthermore, for an integer $k\geq 2$, the $k$-fold tensor products $\mu_{p^n}^{\otimes k}$ and $\mu_{p^\infty}^{\otimes k}$ can be endowed with $\Gal(\bar{F}/F)$-module structure via the diagonal action. For convenience, we shall sometimes write $\mu_{p^\infty}^{\otimes 1}= \mu_{p^\infty}$ and $\mu_{p^\infty}^{\otimes 0} =\Qp/\Zp$, where $\Qp/\Zp$ is understood to have a trivial $\Gal(\bar{F}/F)$-action. Write $H^{i}_{\acute{e}t}\big(F, \mu_{p^n}^{\otimes k}\big) = H^{i}_{\acute{e}t}\big(\mathrm{Spec}(F), \mu_{p^n}^{\otimes k}\big)$ for the \'etale cohomology groups, where  $\mu_{p^n}^{\otimes k}$ is viewed as an \'etale sheaf over the scheme $\mathrm{Spec}(F)$ in the sense of \cite[Chap.\ II]{Mi}. The direct limit of these groups is then denoted by
\[ H^{i}_{\acute{e}t}\big(F, \mu_{p^\infty}^{\otimes k}\big) := \ilim_n H^{i}_{\acute{e}t}\big(F, \mu_{p^n}^{\otimes k}\big). \]

We can now call upon the (Bloch-Lichtenbaum) motivic spectral
sequence
    \[ E^{rs}_2 = H^{r-s}_{\acute{e}t}\big(F, \mu_{p^i}^{\ot (-s)}\big) \Longrightarrow K_{-r-s}(F,\Z/p^i) \]
    (cf. \cite[Chap VI, Theorem 4.2]{WeiKbook}), where $s\leq r\leq 0$. (Strictly speaking, the $E^{rs}_2$-terms should be motivic cohomology groups, which is why the spectral sequence is also coined the ``motivic-to-$K$-theory spectral sequence". Of course, we now know that these motivic cohomology groups are isomorphic to the \'etale cohomology group in this field settings thanks to the work of Rost-Voevodsky \cite{Vo}. Interested readers are referred to \cite{HW} and \cite[Chap.\ VI, Historical Remark 4.4]{WeiKbook}, and the references therein for the history of this monumental spectral sequence.)
    Taking direct limit of the spectral sequence with respect to $i$, we obtain
    \[ H^{r-s}_{\acute{e}t}\big(F, \mu_{p^\infty}^{\ot (-s)}\big) \Longrightarrow K_{-r-s}(F,\Qp/\Zp). \]
   It is a standard fact that $H^{k}_{\acute{e}t}\left(F, -\right) =0$ for $k\geq 3$ (cf.\ \cite[Proposition 8.3.18]{NSW}). For $s\leq -2$, we have $H^{2}_{\acute{e}t}\big(F, \mu_{p^\infty}^{\ot (-s)}\big)=0$ by an observation of Lichtenbaum \cite[Lemma 9.5]{Lic}. Hence the spectral sequence degenerates to two diagonal lines $r-s=0$ and $r-s=1$. This degenerating supplies an identification
   \[ H^1_{\acute{e}t}\big(F, \mu_{p^\infty}^{\ot (n+1)}\big)\cong K_{2n+1}(F,\Qp/\Zp).\]
   One has a similar identification for $L$. These groups are in turn linked via the Tate spectral sequence \cite[Theorem 2.5.3]{NSW}
\[ H_r\left(G, H^{-s}_{\acute{e}t}\big(L, \mu_{p^\infty}^{\ot (n+1)}\big)\right)\Longrightarrow H^{-r-s}_{\acute{e}t}\big(F, \mu_{p^\infty}^{\ot (n+1)}\big). \]
As seen above, $H^{k}_{\acute{e}t}\big(L, \mu_{p^\infty}^{\ot (n+1)}\big)=0$ for $k\geq 2$. Therefore, the initial term of the spectral sequence is at the coordinate $(0,-1)$ and reading off this term, we have
\[ H^{1}_{\acute{e}t}\big(L, \mu_{p^\infty}^{\ot (n+1)}\big)_G\cong H^{1}_{\acute{e}t}\big(F, \mu_{p^\infty}^{\ot (n+1)}\big),\]
which gives the required isomorphism. The proof of the proposition is now complete.
\epf

It follows immediately from Lemma \ref{torsionK2} and Proposition \ref{surjective main at odd p} that
\[\coker\Big(N_{L/F}:K_{2n}(L)\lra K_{2n}(F)\Big) = \coker\Big(N_{L/F}:K_{2n}(L)[2^\infty]\lra K_{2n}(F)[2^\infty]\Big).\]
It therefore remains to analyse the latter which has been accomplished by {\O}stv{\ae}r (see \cite[Theorem 1.2]{Os}), and we shall state his result here.

\bt[{\O}stv{\ae}r] \label{Ostvaer}
Suppose that either $n\neq 1~(mod~4)$ or no real primes of $F$ is ramified in $L$. Then the map $N_{L/F}:K_{2n}(L)[2^\infty]\lra K_{2n}(F)[2^\infty]$ is surjective.

Suppose that $n= 1~(mod~4)$ and some real primes of $F$ is ramified in $L$. Then there
is an exact sequence
\[ K_{2n}(L)[2^\infty]\stackrel{N_{L/F}}{\lra} K_{2n}(F)[2^\infty]\lra  \prod_{v\in S_r}\Z/2 \lra 0,\]
where $S_r$ is the set of real primes of $F$ ramified in $L$.
\et

We are now in position to prove our main results.

\bpf[Proof of Theorems \ref{main combine} and \ref{main cor}]
If $F$ has no real primes ramifying in $L$, then the direct sum is empty. In the event that $n\neq 1$ (mod 4), Suslin's result \cite[Theorem 4.9]{Sus} tells us that $K_{2n}(F_v)_{tor}=0$ for every real prime $v$. Therefore, in either of these situations, the exact sequence in question is saying that $N_{L/F}$ is surjective. But this follows immediately from a combination of Proposition \ref{surjective main at odd p} with Theorem \ref{Ostvaer}.

For the remainder of the proof, we may assume that $n= 1~(\mathrm{mod}~4)$ and some real primes of $F$ are ramified in $L$. We now set to establish Theorem \ref{main cor} first.

Suppose that $x\in K_{2n}(F)$ is a norm from $K_{2n}(L)$.
Let $v$ be any real prime of $F$ that is ramified in $L$ and let $w$ be a complex prime of $L$ which lies above $v$. Consider the following commutative diagram
\[   \entrymodifiers={!! <0pt, .8ex>+} \SelectTips{eu}{}\xymatrix{
     K_{2n}(L)  \ar[d]^{N_{L/F}} \ar[r] & K_{2n}(L_w)
    \ar[d]^{N_{L_w/F_v}} \\
     K_{2n}(F)\ar[r]^{l_v} & K_{2n}(F_v)  } \]
By Suslin's result \cite[Theorem 4.9]{Sus}, the group $K_{2n}(L_w)\cong K_{2n}(\mathbb{C})$ is uniquely divisible and so has no torsion. But as seen in Lemma \ref{torsionK2}, $K_{2n}(L)$ is a torsion group, and so the top horizontal map has to be the trivial homomorphism. Hence it follows that $l_v(x) = 0$.

Conversely, suppose that $l_v(x)=0$ for every real prime $v$ of $F$ is ramified in $L$. Write $x = x_2 + x'$, where $x_2\in K_{2n}(F)[2^\infty]$ and $x'$ an element of odd order. By Proposition \ref{surjective main at odd p}, there exists $y'\in K_{2n}(L)$ such that $N_{L/F}(y') = x'$. It therefore remains to show that $x_2$ has a preimage under $N_{L/F}$. For this, we shall make use of {\O}stv{\ae}r's criterion (see \cite[Theorem 1.1]{Os}). In other words, we need to show that the image of $x_2$ in $K_{2n}(F_u)[2^\infty]$ is in the image of the norm map
\[ \bigoplus_{w|u}K_{2n}(L_w)[2^\infty]\lra K_{2n}(F_u)[2^\infty]\]
for every prime $u$ of $F$. However, as seen in \cite[Claim (2) in Page 492]{Os}, this norm map is surjective for all primes except at the real primes which are ramified in $L$. But for these remaining primes, since we are assuming that their image in $K_{2n}(F_v)[2^\infty]$ is trivial, they are automatically in the image of the norm map. Thus, we may apply {\O}stv{\ae}r's criterion to conclude that $x_2$ comes from the norm of an element in $K_2(L)[2^\infty]$. This finishes the proof of Theorem \ref{main cor}.

We finally come to proving Theorem \ref{main combine} for the remaining case. Observe that Theorem \ref{main cor} is saying that the sequence in question is exact at $K_{2n}(F)$. This in turn implies that the homomorphism $l$ induces an injection
\[  \coker(N_{L/F}) \hookrightarrow\bigoplus_{v\in S_r}K_{2n}(F_v)_{tor}. \]
Now, by virtue of Proposition \ref{surjective main at odd p}, Theorem \ref{Ostvaer} and Suslin's result \cite[Theorem 4.9]{Sus}, these two groups have the same order, and so the injection has to be an isomorphism. Consequently, the homomorphism $l$ is surjective, as required. The proof of Theorem \ref{main combine} is therefore complete.
\epf

\section{Further Remarks}

In view of Theorem \ref{main cor}, it would be of interest to examine the homomorphism $l_v: K_{2n}(F)\lra  K_{2n}(F_v)$ for a real prime $v$ of $F$ in greater depth. In particular, we would like to understand when an element $x\in K_{2n}(F)$ is mapped to zero in $K_{2n}(F_v)$. We describe this for the case of $n=1$. At this point of writing, we do not know how to approach this problem for $n>1$ with $n=1$ (mod 4).

As before, let $F$ be a number field and $v$ a real prime of $F$. For $f,g\in F^\times$, consider the element $\{f,g\}\in K_2(F)$, where $\{-,-\}$ is the Steinberg symbol in $F$. Under the localization map $K_{2}(F)\lra  K_{2}(F_v)$, the symbol $\{f,g\}$ is mapped to $\{l_v(f),l_v(g)\}_v$, where $\{-,-\}_v$ is the Steinberg symbol in $F_v=\mathbb{R}$. As seen previously, the element $\{l_v(f),l_v(g)\}_v$ must lie in $K_2(F_v)_{tor}=\Z/2$. Furthermore, it follows from \cite[Chap.\ III, Example 6.2.1]{WeiKbook} that $\{l_v(f),l_v(g)\}_v \neq 0$ if and only if $l_v(f),l_v(g)\in F_v=\mathbb{R}$
are both negative. We can now record these findings in the form of the following proposition.

\bp \label{normK2}
 Suppose that $L/F$ is a finite extension of number fields. Let $f,g\in F^\times$.  Then $\{f,g\}\in K_2(F)$ is a norm of an element in $K_{2}(L)$ if and only if for every real prime $v$ of $F$ that is ramified in $L$, at least one of $l_v(f)$ and $l_v(g)$ is non-negative.
\ep

The preceding proposition gives a criterion in determining when a symbol can be realized as a norm. In principle, we can build on these ideas to determine whether an arbitrary element $x$ of $K_2(F)$ is a norm element. Indeed, by Matsumoto's theorem \cite[Chap.\ III, Theorem 6.1]{WeiKbook}, the element $x$ can be expressed as a finite product of Steinberg symbols $\{f,g\}$. For each real prime $v$, we have seen above that $l_v(\{f, g\}) \neq 0$ if and only if $l_v(f)$ and $l_v(g)$ are negative. But as long as the number of symbols with both $l_v(f)$ and $l_v(g)$ being negative is even, we still have $l_v(x)=0$. We therefore record our observation as follows.

\bp \label{normK2x}
 Let $L/F$ be a finite extension of number fields. Suppose that $x\in K_2(F)$ has the following representation
 \[ x= \prod_{i\in I}\{f_i, g_i\}, \]
 where $I$ is a finite indexing set. Then $x$ is a norm of an element in $K_{2}(L)$ if and only if for each real prime $v$ of $F$ that ramifies in $L$, the number of symbols with both $l_v(f_i)$ and $l_v(g_i)$ being negative is even.
\ep

We end with some examples to illustrate these propositions. We should mention that the propositions and examples in this section could have been
deduced from the work of Bak-Rehmann (we thank the referees for pointing this out).

\begin{itemize}
           \item[$(1)$] Let $L$ be a number field with at least one pair of complex primes. Then it follows from Proposition \ref{normK2} that $\{-1,-1\}\in K_2(\Q)$ cannot be realized as a norm of an element in $K_2(L)$.

           \item[$(2)$] The symbol $\{3,-1\}\in K_2(\Q)$ can be realized as a norm of an element in $K_2\big(\Q(\sqrt{-1})\big)$. It is interesting to note that the elements $3$ and $-1$ themselves are not norm elements of $\Q(\sqrt{-1})$.

           \item [$(3)$] Let $m$ be a squarefree integer $>1$. A classical result of Bass-Tate (for instance, see \cite[Chap. III, Lemma 6.1.4]{WeiKbook}) asserts that $K_2\big(\Q(\sqrt{m})\big)$ is generated by symbols of the form $\{a, \sqrt{m}-b\}$ and $\{c,d\}$ with $a,c,d\in\Q^\times$ and $b\in\Q$. In the remaining discussion, we shall concentrate on symbols of the form $\{a, \sqrt{m}-b\}$.

               If $a>0$, then $l_v(\{a, \sqrt{m}-b\})=0$ for every real prime $v$ of $\Q(\sqrt{m})$.

               Now suppose that $a<0$. We shall write $l_v$ and $l_{v'}$ for the embedding of $F$ to $\mathbb{R}$ given by $\sqrt{m}\mapsto \sqrt{m}$ and  $\sqrt{m}\mapsto -\sqrt{m}$ respectively. We then consider three cases: (i) $\sqrt{m}>b$, (ii)  $-\sqrt{m} <b <\sqrt{m}$ and (iii) $b < -\sqrt{m}$.  In case (i), one easily checks that $l_v(\{a, \sqrt{m}-b\})=0$ and $l_{v'}(\{a, \sqrt{m}-b\})=0$. In case (ii), we see that $l_v(\{a, \sqrt{m}-b\})=0$ but $l_{v'}(\{a, \sqrt{m}-b\})\neq 0$. Finally, for case (iii), we have $l_v(\{a, \sqrt{m}-b\})\neq0$ and $l_{v'}(\{a, \sqrt{m}-b\})\neq 0$.

                As a further illustration, we let $L$ be an extension of $\Q(\sqrt{m})$ at which both real primes of $\Q(\sqrt{m})$ are ramified (for instance, one may take $L$ to be $\Q(\sqrt{m}, \sqrt{-1})$, $\Q(\sqrt{m}, \zeta_5)$ etc.). Then the above discussion shows that $\{a, \sqrt{m}-b\}$ is a norm element from $K_2(L)$ if and only if either $a>0$ or $\sqrt{m}>b$.

           \item [$(4)$] Let $m$ be a cubefree integer $>1$. Then the field $\Q(\sqrt[3]{m})$ has only one real prime which is the inclusion $\Q(\sqrt[3]{m})\subseteq \mathbb{R}$. Let $L$ be a finite extension of $\Q(\sqrt[3]{m})$ at which this unique real prime is ramified. Then the symbol $\{f,g\}\in K_2\big(\Q(\sqrt[3]{m})\big)$ can be realized as a norm of an element in $K_2(L)$ if and only if either $f>0$ or $g>0$.
               \end{itemize}

\subsection*{Acknowledgement}
The author likes to thank the anonymous referees for the many helpful comments and suggestions.
This research is supported by the
National Natural Science Foundation of China under Grant No. 11771164 and the Fundamental Research Funds for the Central Universities of CCNU
under grant CCNU20TD002.


\footnotesize
\begin{thebibliography}{00}

\bibitem{BR} A. Bak and U. Rehmann, $K_2$-analogs of Hasse's norm theorems.
Comment. Math. Helv. 59 (1984), no. 1, 1-11.

\bibitem{Bo} A. Borel, Stable real cohomology of arithmetic groups. Ann. Sci. \'Ecole Norm. Sup. (4) 7 (1974), 235-272.

\bibitem{Br} W. Browder,  Algebraic $K$-theory with coefficients $\Z/p$. Geometric applications of homotopy theory I, pp. 40-84,
Lecture Notes in Math., 657, Springer, Berlin, 1978.

\bibitem{CT} J.-L. Colliot-Th\'el\`{e}ne, Hilbert's Theorem 90 for $K_2$, with application to the Chow groups of rational surfaces. Invent. Math. 71 (1983), no. 1, 1-20.

\bibitem{Gar} H, Garland, A finiteness theorem for $K\sb{2}$ of a number field. Ann. of Math. (2) 94 (1971), 534--548.

\bibitem{HW} C. Haesemeyer and C. Weibel, The norm residue theorem in motivic cohomology.
Annals of Mathematics Studies, 200. Princeton University Press, Princeton, NJ, 2019.



\bibitem{Lic} S. Lichtenbaum, On the values of zeta and $L$-functions. I.
Ann. of Math. (2) 96 (1972), 338-360.

\bibitem{MS} A. S. Merkur'ev and A. A. Suslin, $K$-cohomology of Severi-Brauer varieties and the norm residue homomorphism. (Russian) Izv. Akad. Nauk SSSR Ser. Mat. 46 (1982), no. 5, 1011-1046.

\bibitem{Mi} J. Milne, Arithmetic Duality Theorems. Second edition. BookSurge, LLC, Charleston, SC, 2006. viii+339 pp.


\bibitem{NSW} J. Neukirch; A. Schmidt; K. Wingberg,
Cohomology of Number Fields, 2nd edn., Grundlehren Math.
Wiss. 323 (Springer-Verlag, Berlin, 2008).

\bibitem{OVV} D. Orlov; A. Vishik;  V. Voevodsky, An exact sequence for $K^M_*/2$ with applications to quadratic forms, Ann. of Math. (2) 165 (2007), no. 1, 1-13.

\bibitem{Os}  P. A. {\O}stv{\ae}r, A norm principle in higher algebraic $K$-theory.
Bull. London Math. Soc. 35 (2003), no. 4, 491-498.

\bibitem{Sou} C. Soul\'e, $K$-th\'eorie des anneaux d'entiers de corps de nombres et cohomologie \'etale. Invent. Math. 55 (1979), no. 3, 251-295.

\bibitem{Sus} A. A. Suslin, On the $K$-theory of local fields. Proceedings of the Luminy conference on algebraic $K$-theory (Luminy, 1983). J. Pure Appl. Algebra 34 (1984), no. 2-3, 301-318.

\bibitem{Qui72} D. Quillen, On the cohomology and $K$-theory of the general linear groups over a finite field, Ann. of Math. (2) 96 (1972), 552-586.

\bibitem{Qui73a} D. Quillen, Higher algebraic $K$-theory. I. Algebraic $K$-theory, I: Higher $K$-theories (Proc. Conf., Battelle Memorial Inst., Seattle, Wash., 1972), pp. 85-147. Lecture Notes in Math., Vol. 341, Springer, Berlin 1973.

\bibitem{Qui73b} D. Quillen, Finite generation of the groups $K\sb{i}$ of rings of algebraic integers. Algebraic $K$-theory, I: Higher $K$-theories (Proc. Conf., Battelle Memorial Inst., Seattle, Wash., 1972), pp. 179-198. Lecture Notes in Math., Vol. 341, Springer, Berlin, 1973.


\bibitem{Vo} V. Voevodsky, On motivic cohomology with $\mathbf Z/l$-coefficients. Ann. of Math. (2) 174 (2011), no. 1, 401-438.


\bibitem{WeiKbook} C. Weibel, The $K$-book. An introduction to algebraic $K$-theory. Graduate Studies in Mathematics, 145. American Mathematical Society, Providence, RI, 2013. xii+618 pp.

\end{thebibliography}
\end{document}